\documentclass[preprint,twocolumn,5p,11pt,authoryear]{elsarticle}

\usepackage{amssymb}
\usepackage{color}
\usepackage{makecell}
\usepackage{multirow}
\usepackage{xurl}

\journal{Journal of Process Control}
\usepackage{algpseudocode}
\usepackage{algorithm}
\usepackage{booktabs}
\usepackage{amsmath}
\usepackage{bm}
\newcommand{\rmt}{\mathrm{T}}
\usepackage{algpseudocode}

\newcounter{remarkcount}
\newcounter{theoremcount}
\newcounter{assumptioncount}
\usepackage{float}
\usepackage[caption = false]{subfig}

\newcommand{\remark}[1]{%
    \stepcounter{remarkcount}
    \textbf{Remark \arabic{remarkcount}:} #1
}

\newcommand{\assumption}[1]{%
    \stepcounter{assumptioncount}
    \textit{Assumption \arabic{assumptioncount}:} #1
}

\newcommand{\theorem}[1]{%
    \stepcounter{theoremcount}
    \textbf{Theorem \arabic{theoremcount}:} #1
}

\begin{document}

\begin{frontmatter}



\title{On Data-Driven Robust Optimization With Multiple Uncertainty Subsets:\\ Unified Uncertainty Set Representation and Mitigating Conservatism}

\author[a]{Yun Li\corref{cor1}}
\ead{y.li-39@tudelft.nl}
\cortext[cor1]{Corresponding author.}

\author[b]{Neil Yorke-Smith}
\author[a]{Tamas Keviczky}

\address[a]{Delft Center for Systems and Control, 
            Delft University of Technology, 
            Delft,
            the Netherlands}
\address[b]{STAR Lab, Delft University of Technology, Delft, the Netherlands}

\begin{abstract}
Constructing uncertainty sets as unions of multiple subsets has emerged as an effective approach for creating compact and flexible uncertainty representations in data-driven robust optimization (RO). This paper focuses on two separate research questions. The first concerns the computational challenge in applying these uncertainty sets in RO-based predictive control. To address this, a monolithic mixed-integer representation of the uncertainty set is proposed to uniformly describe the union of multiple subsets, enabling the computation of the worst-case uncertainty scenario across all subsets within a single mixed-integer linear programming (MILP) problem. The second research question focuses on mitigating the conservatism of conventional RO formulations by leveraging the structure of the uncertainty set. To achieve this, a novel objective function is proposed to exploit the uncertainty set structure and integrate the existing RO and distributionally robust optimization (DRO) formulations, yielding less conservative solutions than conventional RO formulations while avoiding the high-dimensional continuous uncertainty distributions and incurring high computational burden typically associated with existing DRO formulations. Given the proposed formulations, numerically efficient computation methods based on \textit{column-and-constraint generation} (CCG) are also developed. Extensive simulations across three case studies are performed to demonstrate the effectiveness of the proposed schemes.

\end{abstract}


\end{frontmatter}


\section{Introduction}
Decision-making under uncertainty has gained increasing attention in both scientific research and engineering applications. One of the solutions for handling uncertainties in decisions is robust optimization (RO). By representing uncertainties via uncertainty sets, RO can ensure constraint satisfaction for all possible uncertainties without incurring an unaffordable computational burden and requiring the complete information of uncertainty distributions \citep{ben1998robust,ben2009robust, el1998robust}. Despite the popularity and effectiveness of RO formulations, it is well-known that the optimal solution of RO problems can be conservative owning to two main factors. 

The first factor contributing to the conservatism of the robust optimal solution is how the uncertainty set is represented. In recent years, with the availability of abundant measurement data and machine learning techniques, data-driven approaches have been developed to construct compact uncertainty sets by uncovering and exploiting the latent features of the uncertainty data, see \cite{shang2017data,shang2019data,ning2017data,ning2018data,bertsimas2018data,cheramin2021data} and references therein. Among these data-driven approaches, one innovative methodology is constructing the uncertainty set as a union of multiple subsets, see \cite{ning2017data,ning2018data,li2024machine} and references therein. By adopting such a methodology, the data-driven uncertainty set is able to adapt to irregular uncertainty distributions and ensure compactness, which is beneficial for reducing the conservatism of the corresponding robust optimal solutions. 

It should be noted that, while the methodology of representing uncertainty sets as unions of multiple subsets has been proven to be effective in several applications, such as power system operation \citep{zhao2022sustainable} and building climate control \citep{hu2023multi}, this uncertainty set representation together with the existing \textit{column-and-constraint generation} (CCG)-based algorithm could result in computationally demanding optimization problems when applied in robust predictive control settings, which is generally neglected in the existing works. 

One widely adopted algorithmic solution for solving two-stage RO problems is the so-called CCG algorithm \citep{zhao2012exact,zeng2013solving}. When applying the existing CCG algorithm to solve the RO problem with multiple uncertainty subsets, each uncertainty subset entails solving one optimization problem to compute the worst-case uncertainty in each iteration of the CCG algorithm \citep{ning2017data}. As elaborated in Section \ref{sec:predictive_control}, the total number of uncertainty subsets in RO-based predictive control problems could increase exponentially with the length of the prediction horizon, which leads to computationally demanding or even intractable problems and limits the applicability of the existing approaches. 

The second major factor influencing the conservatism of RO solutions is how the objective function is designed. In conventional RO formulation, the objective function represents the worst-case performance, which can lead to conservative optimal solutions \citep{ben1998robust,ben2009robust,bertsimas2022robust,ning2017data}. For mitigating the conservatism of the optimal solution, some alternative objective functions are proposed, such as minimizing the expected performance in stochastic optimization (SO) \citep{birge2011introduction,calafiore2006scenario,campi2008exact} and optimizing the worst-case expected performance in distributionally robust optimization (DRO) \citep{hanasusanto2015distributionally,mohajerin2018data,bertsimas2010models,bertsimas2022robust,wei2020tutorials,kuhn2024dro}. 

However, despite the effectiveness of these existing objective functions, they either require approximating the probability distribution of uncertainties, which could be high-dimensional and is non-trivial to obtain in practice, or entail solving a computationally demanding optimization problem. For example, in Moment-based DRO formulation \citep{hanasusanto2015distributionally}, a semi-definite programming problem needs to be solved; in scenario-based SO \citep{calafiore2006scenario} and Wasserstein-based DRO formulation \citep{mohajerin2018data}, the consideration of a large number of uncertainty scenarios leads to large-scale optimization problems and increased computational burden. Furthermore, the existing design of conservatism-reduced objective functions does not exploit the structure of uncertainty sets, especially the promising option of representing the uncertainty set as a union of multiple subsets in data-driven RO settings.

By considering that the uncertainty set consists of multiple subsets, this paper mainly investigates two separate research questions: 1) \textit{how to properly represent the uncertainty set to deal with the computational challenge when applying it in RO-based predictive control problems}; and 2) \textit{how to exploit the uncertainty set structure to design new objective function for mitigating the conservatism of the robust optimal solution}. The main contributions of this paper are summarized as follows:
\begin{itemize}
    \item A monolithic representation of the uncertainty set consisting of multiple subsets is proposed for RO-based predictive control problems. Given the proposed representation, a computationally efficient CCG algorithm is developed such that only a single MILP problem needs to be solved to compute the worst-case uncertainty scenario across all uncertainty subsets regardless of the length of the prediction horizon. In comparison, conventional solutions entail solving a set of optimization problems whose number grows exponentially with the prediction horizon \citep{ning2017data,zeng2013solving}. 
    \item Aiming at reducing the conservatism for optimizing the worst-case performance in conventional RO formulations, a novel objective function exploiting the structure of the uncertainty set together with a CCG-based computation method are proposed. By combining conventional RO and DRO formulations, the proposed objective function not only achieves less conservative solutions than conventional RO formulations but also avoids approximating the continuous and possibly high-dimensional uncertainty distribution and incurring high computational burden as with the existing SO and DRO formulations, especially in the case of fully adaptive recourse decision and large numbers of uncertainty samples.
    \item Numerical experiments of three case studies are performed to extensively illustrate the issues of the existing approaches with the uncertainty set consisting of multiple subsets and demonstrate the effectiveness of the proposed schemes. 
\end{itemize}

The remaining parts of this paper are organized as follows. Section \ref{sec:preliminary} briefly describes the preliminaries and research gaps. Section \ref{sec:predictive_control} clarifies the issues of the existing approach and proposes a new uncertainty set representation with a CCG-based computation method. In Section \ref{sec:stochastic_dro}, a novel objective function exploiting the structure of the uncertainty set is proposed to mitigate the conservatism of robust optimal solutions. Given the proposed objective function, a CCG-based computation method is also derived. Section \ref{sec:simulation} presents numerical experiment results. Section \ref{sec:conclusion} concludes this paper.

\textbf{Notation}: boldface capital letters denote matrices, and boldface lowercase letters denote vectors. Sets are represented by calligraphic capital letters. $\mathbf{e}$ denotes all-ones vectors with appropriate dimensions. $\otimes$ is the Cartesian product. $\mathbb{E}[\cdot]$ indicates the expectation of a given random variable. 

\section{Preliminaries}\label{sec:preliminary}
This section introduces some preliminaries of the two-stage linear RO problem to be investigated, and briefly discusses the research gaps when considering an uncertainty set consisting of multiple subsets. 

The following two-stage linear RO problem is considered in this work
\begin{subequations}\label{eq:ro}
    \begin{align}
    \min_{\bf x}\ & \mathbf{c}^\rmt\mathbf{x} + \max_{\mathbf{v}\in\mathcal{V}}\min_{\mathbf y}\mathbf{b}^\rmt \mathbf{y} \label{eq:worst_obj}\\
    \text{s.t. } & \mathbf{Ax \leq q},\label{eq:1b}\\
    & \mathbf{Tx + Wy + Mv \leq h},\label{eq:linear_cons}
    \end{align}
\end{subequations}
where $\mathbf{x}\in\mathbb{R}^{p}$ are the first-stage decision variables (interchangeably called \textit{here-and-now} decision variables), which can contain both continuous and binary elements, $\mathbf{y}\in\mathbb{R}^q$ are the second-stage decision variables (interchangeably called \textit{wait-and-see} decision variables or \textit{recourse} decision variables), which are assumed to be continuous, $\mathbf{v}\in\mathbb{R}^m$ denotes the uncertainties with $\mathcal{V}$ as the corresponding uncertainty set, $(\mathbf{c},\mathbf{b},\mathbf{q},\mathbf{h})$ and $(\mathbf{A},\mathbf{T},\mathbf{W},\mathbf{M})$ are parameter vectors and matrices, respectively, with appropriate dimensions. Without loss of generality, in the remaining parts of this paper, the above two-stage RO problem is assumed to be relatively complete recourse, i.e., the optimization problem is feasible for any possible $\mathbf{x}$ and $\mathbf{v}$ \citep{zeng2013solving}. For problems that do not satisfy this assumption, some other strategic solutions discussed in \cite{zhao2012exact} are also applicable in our proposed design. 

To solve the above two-stage RO problem, one commonly adopted approach is the so-called \textit{column-and-constraint generation} (CCG) algorithm. In CCG algorithm, the optimal solution of the RO problem is computed by iteratively solving a master problem and several subproblems (the number of the subproblems to be solved is dependent on the number of uncertainty subsets) to update the lower bound and upper bound of the objective function until a predefined optimality gap is achieved \citep{zeng2013solving,ning2017data}. In this paper, unless mentioned specifically, the existing approach for solving RO problem \eqref{eq:ro} refers to the CCG algorithm.

For reducing the conservatism of the robust optimal solutions, including but not limited to the two-stage RO problem in \eqref{eq:ro}, data-driven approaches exploiting historical uncertainty data and machine learning techniques have been extensively explored for constructing compact uncertainty set. Among these data-driven approaches, one promising methodology is representing the uncertainty set as a union of multiple subsets, such as polytopes and ellipsoids, which has been shown to be effective in improving the compactness and flexibility of the uncertainty sets in handling irregular and complex uncertainty distributions. An uncertainty set $\mathcal{V}$ represented as a union of multiple subsets can be expressed as
\begin{equation}\label{eq:union_set}
    \mathcal{V}:=\bigcup\nolimits_{k=1}^K\mathcal{V}_k
\end{equation}
where $\mathcal{V}_k$ denotes the $k$-th uncertainty subset with $K$ as the total number of subsets. The uncertainty subsets can be any basic sets such as box, ellipsoid, and polytope.

While the existing works have demonstrated the effectiveness of formulating an uncertainty set as a union of multiple subsets as in \eqref{eq:union_set}, there are two research gaps accompanying with such uncertainty set representation to be filled.
\begin{itemize}
    \item \textit{RO-based Predictive Control With Multiple Uncertainty Subsets}: 
    The RO formulation \eqref{eq:ro} has been applied in predictive control problems due to its effectiveness in ensuring robust constraint satisfaction \citep{shang2019data}. However, when employing the data-driven uncertainty set representation \eqref{eq:union_set}, the number of uncertainty subsets could grow exponentially with the prediction horizon. This exponential growth leads to a rapid increase in the number of subproblems solved in each iteration of the existing CCG-based algorithms, making the optimization problem computationally demanding or even intractable. Therefore, developing computationally efficient solutions to address this challenge is crucial for enhancing the practicality of the data-driven RO design with the uncertainty set representation \eqref{eq:union_set} in robust predictive control.
    \item \textit{Conservatism-Reducted Objective Function With Multiple Uncertainty Subsets}: The conventional RO formulation could suffer from conservatism due to the optimization of worst-case performance, as shown in \eqref{eq:ro}. Given the existing RO formulations with multiple uncertainty subsets, the structure of the uncertainty set is not fully considered in designing the objective function. Consequently, how to properly exploit the property that the uncertainty set consists of multiple subsets to design novel RO formulations is beneficial to reduce the conservatism of the robust optimal solution.
\end{itemize}
The above research gaps together with our proposed solutions will be further elaborated in detail in the subsequent sections.

\section{RO-based Predictive Control With Multiple Uncertainty Subsets}\label{sec:predictive_control}
This section details the research gap when applying the uncertainty set structure \eqref{eq:union_set} in robust predictive control, and proposes a novel mixed-integer representation to uniformly describe the unified uncertainty set such that the computational efficiency of the CCG algorithm for solving the RO problem is remarkably improved. 

\subsection{Problem Formulation}
Theoretically, even with the uncertainty set represented as in \eqref{eq:union_set}, existing data-driven methods for constructing uncertainty sets and the solutions for solving \eqref{eq:ro} can be directly extended to robust predictive control problems. However, practical challenges arise when applying these methods, which have been overlooked in existing works. 

For simplicity of illustration, consider the following deterministic formulation of a predictive control problem
\begin{subequations}\label{eq:predictive_control}
    \begin{align}
        \min_{\mathbf{u}_t}\ &\sum_{t=1}^{N}l(\mathbf{s}_t,\mathbf{u}_t) \\
        \text{s.t. } & \mathbf{s}_{t+1} = \bm{\Phi}\mathbf{s}_t + \bm{\Gamma}\mathbf{u}_t + \mathbf{v}_t\\
        & \mathbf{s}_t\in\mathcal{S}, \mathbf{u}_t\in\mathcal{U},\ \mathbf{v}_t\in\mathcal{V}_t\\
        & t = 1,\cdots,N
    \end{align}
\end{subequations}
where $\mathbf{s}_t$ denotes the system states vector, $\mathbf{u}_t$ denotes the control input vector, $\mathbf{v}_t$ the uncertainties, $l(\mathbf{s}_t,\mathbf{u}_t)$ the stage cost function, $\mathcal{S}$ and $\mathcal{U}$ are feasible sets of system states and control inputs, respectively, $\mathcal{V}_t$ is the uncertainty set of $\mathbf{v}_t$, subscript $t$ denotes the $t$-th time step, and $N$ the length of the prediction horizon. In predictive control settings, it is common that the uncertainties $\mathbf{v}_t\ (t=1,\cdots,N)$ are assumed to be independent and identically distributed (I.I.D.) \citep{calafiore2012robust,mesbah2016stochastic}.

After reformulating the predictive control problem \eqref{eq:predictive_control} into the corresponding RO formulation \eqref{eq:ro}, the uncertainty vector $\mathbf{v}$ is the stacked uncertainty sequences $\mathbf{v}_t$ within the prediction horizon, i.e., $\mathbf{v}:=[\mathbf{v}_1^\rmt,\cdots,\mathbf{v}_N^\rmt]^\rmt$. In the existing data-driven RO framework, the uncertainty set $\mathcal{V}$ is directly constructed for the stacked uncertainty vector $\mathbf{v}$. While this approach is still theoretically feasible, it will incur some issues. On one hand, for a fixed number of uncertainty samples $\mathbf{v}_t$, considering the stacked uncertainty $\mathbf{v}$ reduces the size of the training dataset (shrinking it to $\frac{1}{N}$ of the original dataset for $\mathbf{v}_t$). This reduction may result in insufficient data and degrade the performance of data-driven approaches for constructing uncertainty sets. On the other hand, directly modelling the uncertainty set for the stacked uncertainty $\mathbf{v}$ increases the dimensionality of uncertainty samples, particularly for long prediction horizon $N$. As shown in existing studies, high-dimensional uncertainties pose challenges in implementing data-driven methods, including difficulties in hyper-parameter tuning and performance validation.

For example, in building climate predictive control, the prediction horizon could span 12 hours with a sampling period of $30$ minutes, and a typical uncertainty is the prediction error of ambient temperature. Given one year of historical uncertainty data of $\mathbf{v}_t$, the number of uncertainty samples for $\mathbf{v}_t$ is 365*24*2. However, when considering the stacked uncertainty samples for $\mathbf{v}=[\mathbf{v}_t,\cdots,\mathbf{v}_N]^\rmt$ with $N= 24$, the training set size reduces to $365*2$. Consequently, the dataset for $\mathbf{v}$ may be too small to ensure the effectiveness of data-driven approaches in constructing uncertainty sets.

Under the I.I.D. assumption of $\mathbf{v}_t$, one solution is to construct a data-driven uncertainty set for $\mathbf{v}_t$, and then extend it to $\mathbf{v}$. Denoting the $k$-th uncertainty subsets of $\mathbf{v}_t$ as $\mathcal{V}_{t,k}$, the uncertainty set for $\mathbf{v}=[\mathbf{v}_1^\rmt, \cdots,\mathbf{v}_N^\rmt]^\rmt$, denoted as $\mathcal{V}$, can be expressed as
\begin{equation}\label{eq:uncerset_exp}
        \mathcal{V}= \underbrace{\mathcal{V}_t\otimes\ldots\otimes\mathcal{V}_t}_{N},\quad \mathcal{V}_t = \bigcup\nolimits_{k=1}^K\mathcal{V}_{t,k}
\end{equation}
where $K$ is the total number of uncertainty subsets for $\mathbf{v}_t$, and $N$ is the length of prediction horizon.


\remark{It is evident from \eqref{eq:uncerset_exp} that the number of subsets constituting $\mathcal{V}$ is $K^N$, which grows exponentially with the prediction horizon $N$. By adopting the existing CCG algorithm to solve the RO problem \eqref{eq:ro} with the uncertainty set \eqref{eq:uncerset_exp}, each uncertainty subset requires solving a bilinear or mixed-integer linear optimization problem to compute the worst-case uncertainty in every iteration of the algorithm \citep{ning2017data}. Consequently, handling $K^N$ subsets translates to solving $K^N$ optimization problems per iteration, making the approach computationally demanding or even intractable even for moderate values of $K$ and $N$

Therefore, the computational challenges arising from representing the uncertainty set as a union of multiple subsets in \eqref{eq:uncerset_exp} highlight the need for an effective formulation of the uncertainty set for $\mathbf{v}$. Designing such a formulation is crucial for ensuring that the corresponding RO-based predictive control problem remains computationally tractable and practically applicable.
}
\subsection{Unified Uncertainty Set Representation and Computational Method Design}
This section aims to propose a novel representation to uniformly represent the uncertainty set of $\mathbf{v}$ based on the uncertainty subsets of $\mathbf{v}_t$ and derive a computationally efficient solution for solving the corresponding RO problem. 

\assumption{The uncertainty set of $\mathbf{v}_t$, denoted as $\mathcal{V}_t$, consists of $K$ subsets $\mathcal{V}_{t,k}$ with each uncertainty subset as a nonempty and bounded polytope defined as $\mathcal{V}_{t,k} := \{\mathbf{v}_t|\mathbf{D}_{k}\mathbf{v}_t\leq \mathbf{d}_k\}$. }

On the basis of Assumption 1, the following monolithic mixed-integer formulation of uncertainty set $\mathcal{V}$ is derived
\begin{equation}\label{eq:MI_uncertset}
    \begin{aligned}
    \mathcal{V} := &\bigg\{[\mathbf{v}_1^\rmt,\cdots,\mathbf{v}_N^\rmt]^\rmt\bigg\vert\sum_{k=1}^K\delta_{t,k}\mathbf{D}_k\mathbf{v}_t\leq \sum_{k=1}^K\delta_{t,k}\mathbf{d}_{k},\\&\delta_{t,k}\in\mathbb{B},\ \sum\nolimits_{k=1}^K\delta_{t,k} = 1,\ t = 1,\cdots,N\bigg\}
    \end{aligned}
\end{equation}
where $\delta_{t,k}\in\{0,1\}$ are auxiliary variables. For the formulation \eqref{eq:MI_uncertset}, given any subset of $\mathcal{V}$, a set of feasible variables $\{\delta^*_{t,k}\mid t= 1,\cdots,N, k = 1,\cdots,K\}$ can be found to uniquely represent a specific uncertainty subset. Conversely, any feasible $\{\delta^*_{t,k}\mid t= 1,\cdots,N, k = 1,\cdots,K\}$ will also define an admissible uncertainty subset of $\mathcal{V}$. Instead of explicitly representing the uncertainty set $\mathcal{V}$ as a union of multiple subsets as in \eqref{eq:uncerset_exp}, the proposed mixed-integer representation \eqref{eq:MI_uncertset} gives a monolithic implicit representation of $\mathcal{V}$ that is feasible for any sizes of uncertainty subsets $K$ and prediction horizon $N$.

\theorem{\textit{Considering the uncertainty set \eqref{eq:MI_uncertset} and decomposing the parameter matrix $\mathbf{M}$ into $N$ column blocks $\{\mathbf{M}_t\}_{t=1}^N$ with appropriate dimensions, the RO problem in \eqref{eq:ro} can be reformulated as in \eqref{eq:RO_reform}, which can be solved via Algorithm \ref{alg:ccg_worst_case} within finite iterations.}}
\begin{subequations}\label{eq:RO_reform}
    \begin{align}
        \min_{\mathbf{x}}\ &\mathbf{c}^\rmt\mathbf{x} + \max_{\delta_{t,k},\mathbf{v}_t}\min_{\mathbf{y}}\mathbf{b}^\rmt\mathbf{y}  \\
        \text{s.t. }& \mathbf{A} \mathbf{x} \leq \mathbf{q},\\
        &\mathbf{T}\mathbf{x} + \mathbf{Wy} + \sum\nolimits_{t=1}^{N}\mathbf{M}_t\mathbf{v}_t\leq \mathbf{h},\\
        &\sum\nolimits_{k=1}^K\delta_{t,k}\mathbf{D}_k\mathbf{v}_t\leq \sum\nolimits_{k=1}^K\delta_{t,k}\mathbf{d}_k,\\
        &\delta_{t,k}\in\mathbb{B},\ k=1,\cdots,K,\ t = 1,\cdots,N,\\
        & \sum\nolimits_{k=1}^K\delta_{t,k} = 1,\ t = 1,\cdots,N.
    \end{align}
\end{subequations}
\textit{Proof}: the proof of this theorem is presented in the Appendix.\hfill$\square$

The master problem \textbf{MP} and the subproblem \textbf{SP}1 entailed in Algorithm \ref{alg:ccg_worst_case} are defined as:
\begin{subequations}\label{eq:mp2}
   \begin{align}
    \textbf{MP}:\quad    \min_{\mathbf{x},\eta,\mathbf{y}_i}\ & \mathbf{c}^\rmt \mathbf{x} + \eta \\
    \text{s.t. } & \mathbf{b^\rmt y}_i \leq \eta, \\
    & \mathbf{Ax\leq q},\\
    & \mathbf{Tx +Wy}_i + \sum\nolimits_{t=1}^N\mathbf{M}_t\mathbf{v}_{t} \leq \mathbf{h},\\
    & i = 1,2,\cdots, r,\ k = 1,\cdots,K.
    \end{align}
\end{subequations}
\begin{subequations}\label{eq:subp_reform1}
    \begin{align}
    \textbf{SP}1:&\notag\\
    \quad \max_{\delta_{t,k},\mathbf{v}_t}\min_{\bf y}&\ \mathbf{b}^\rmt \mathbf{y}\\
    \text{s.t. } & \mathbf{Tx+Wy}+\sum\nolimits_{t=1}^N\mathbf{M}_t\mathbf{v}_t \leq \mathbf{h}, \\
    & \sum\nolimits_{k=1}^K\delta_{t,k}\mathbf{D}_k\mathbf{v}_t \leq \sum\nolimits_{k=1}^K\delta_{t,k}\mathbf{d}_k,\label{eq:MI_uncer2}\\
    & \sum\nolimits_{k=1}^K\delta_{t,k} = 1,\ \delta_{t,k}\in\mathbb{B},\\
    & k = 1,\cdots,K,\ t = 1,\cdots,N.
    \end{align}    
\end{subequations}

\begin{algorithm}[htb]
\caption{\textit{column-and-constraint generation} algorithm for solving \eqref{eq:RO_reform}.}
\label{alg:ccg_worst_case}
    \begin{algorithmic}[1]
        \Statex \textbf{Input}: suboptimality gap $\epsilon$
        \Statex \textbf{Output}: the optimal decision variable $\mathbf{x}$ and objective function value $\mathbf{c}^\rmt\mathbf{x}^* + \eta^*$
        \State Set $LB = -\infty$, $UB = \infty$, $r = 0$ 
        \While{$|UB-LB|>\epsilon$}
        \State Solve \textbf{MP} in \eqref{eq:mp2} to derive solutions $(\mathbf{x}^*, \eta^*$) and update $LB = \mathbf{c^\rmt} \mathbf{x}^*+\eta^*$
        \State Solve \textbf{SP}1 in \eqref{eq:subp_reform1} or \textbf{SP}2 in \eqref{eq:subp_reform2} to derive solutions $\{\mathbf{v}^*_{t},\mathbf{y}^*\}$, and update $UB$ as 
        \begin{equation*}
        \begin{aligned}
            UB = \min & \bigg\{UB, \mathbf{c^\rmt} \mathbf{x}^* + \mathbf{b^\rmt}\mathbf{ y}^*\bigg\}
            \end{aligned}
        \end{equation*}
        \State Create decision variables $\mathbf{y}_r$, set parameters $\mathbf{v}_{t,r} = \mathbf{v}_t^*\ (t=1,\cdots,N)$ , and add the following constraints to \textbf{MP} in \eqref{eq:mp2} 
        \begin{align*}
        \begin{cases}
            \mathbf{b}^\rmt \mathbf{y}_t \leq \eta, \\
            \mathbf{Tx + Wy}_r +\sum_t\mathbf{M}_t\mathbf{v}^*_{t,r} \leq \mathbf{h}.
            \end{cases}
        \end{align*}
        
        \State $r \leftarrow r+1$
        \EndWhile
        \State \textbf{Return}: $\mathbf{x}^*$ and $\mathbf{c^\rmt}\mathbf{x}^* + \eta^*$  
    \end{algorithmic}
\end{algorithm}

The subproblem \eqref{eq:subp_reform1} contains bilinear terms $\delta_{t,k}\mathbf{v}_t$ due to the auxiliary binary decision variable $\delta_{t,k}$, and hence is non-convex and might be intractable to some solvers. These bilinear terms can be further reformulated as mixed-integer linear constraints via big-M formulation. The reformulated subproblem \textbf{SP}2 is given in \eqref{eq:subp_reform2}.
\begin{subequations}\label{eq:subp_reform2}
    \begin{align}
    \textbf{SP}2:\notag\\
    \max_{\delta_{t,k},\mathbf{v}_t,\mathbf{w}_{t,k}}&\min_{\mathbf{y}}\ \mathbf{b^\rmt y} \\
        \text{s.t. }& \mathbf{Tx + Wy} + \sum\nolimits_{t=1}^N\mathbf{M}_t\mathbf{v}_t \leq \mathbf{h}, \\
        &\sum\nolimits_{k=1}^K\mathbf{D}_k \mathbf{w}_{t,k} \leq \sum\nolimits_{k=1}^K\delta_{t,k}\mathbf{d}_k, \\
        &-\Delta(1-\delta_{t,k}) \mathbf{e} \leq \mathbf{v}_t - \mathbf{w}_{t,k} \leq \Delta(1-\delta_{t,k})\mathbf{e},\\
        &-\Delta\delta_{t,k}\mathbf{e} \leq \mathbf{w}_{t,k} \leq \Delta\delta_{t,k}\mathbf{e},\\
        &k=1,\cdots,K,\  t= 1,2,\cdots,N
    \end{align}
\end{subequations}
where $\Delta > 0$ is a sufficiently large constant, and $\mathbf{e}$ is an all-ones vector with appropriate dimension.

\remark{In contrast to the existing solution for solving the RO problem \eqref{eq:ro} with uncertainty set \eqref{eq:uncerset_exp}, in which $K^N$ numbers of subproblems need to be solved in each algorithm iteration, the reformulated uncertainty set \eqref{eq:MI_uncertset} enables a monolithic description of the uncertainty set and ensures that only a single subproblem (\textbf{SP}1 or \textbf{SP}2), whose decision variables and constraints increase linearly with $K$ and $N$, to be solved in Algorithm \ref{alg:ccg_worst_case} for computing the worst-case uncertainty. This can remarkably improve the computational efficiency of the proposed approach in case of large values of $K$ and/or $N$. Additionally, the proposed approach enables constructing data-driven uncertainty sets for $\mathbf{v}_t$ instead of the stacked uncertainty $\mathbf{v}$ without incurring a high computational burden for solving the corresponding RO problem \eqref{eq:ro}. In comparison with $\mathbf{v}$, $\mathbf{v}_t$ has a lower dimension and more data samples, which will further boost the performance of the data-driven approaches for constructing uncertainty sets. For example, a lower uncertainty dimension makes it easier to tune the parameters/hyper-parameters and evaluate the performance of the data-driven approaches. A larger dataset can also prevent the data-driven approaches from suffering performance degradation, such as overfitting, underfitting, poor model calibration, etc., due to insufficient training data.}

\section{Mitigating Conservatism via Objective Function Design}\label{sec:stochastic_dro}
This section proposes an alternative objective function for RO problems by leveraging the structure of uncertainty sets composed of multiple subsets. Additionally, an algorithmic solution is developed to efficiently solve the corresponding RO problem.

\subsection{Problem Formulation}
Based on the structure of the uncertainty set -- a union of multiple subsets -- this section presents a new formulation of the objective function for RO problems to reduce the conservatism of the optimal solution. 

The conventional RO formulation in \eqref{eq:ro} could lead to conservative solutions due to the objective of worst-case performance optimization. Many existing works have been done to propose new objective functions to reduce the conservatism of the optimal solution. Two prominent alternative objective functions are given as follows.
\begin{equation}
    \min_{\mathbf{x}}\left\{\mathbf{c}^\rmt\mathbf{x} + Q(\mathbf{v},\mathbf{x}) \right\}
\end{equation}
where two alternative options of $Q(\mathbf{v},\mathbf{x})$ are
\begin{subequations}\label{eq:obj_options}
\begin{align}
    &\text{Option 1: }Q(\mathbf{v},\mathbf{x}) := \mathbb{E}_{f(\mathbf{v})}\left[\min_{\mathbf{y}}\mathbf{b}^\rmt \mathbf{y}\right]\label{eq:Q_mean},\\
    &\text{Option 2: }Q(\mathbf{v},\mathbf{x}) := \max_{f(\mathbf{v})\in\mathcal{P}}\mathbb{E}_{f(\mathbf{v})}\left[\min_{\mathbf{y}}\mathbf{b}^\rmt \mathbf{y}\right]\label{eq:Q_dro}
\end{align}
\end{subequations}
with $f(\mathbf{v})$ as the probability distribution of the uncertainty $\mathbf{v}$, and $\mathcal{P}$ denotes the ambiguity set defining the admissible set of $f(\mathbf{v})$.

The definition in \eqref{eq:Q_mean} optimizes the expected performance w.r.to the probability of uncertainty $\mathbf{v}$, which is generally adopted in the setting of stochastic optimization (SO). By further considering the uncertainty of the probability distribution $f(\mathbf{v})$, the definition \eqref{eq:Q_dro} formulates a so-called DRO problem, which optimizes the worst-case expected performance w.r.to the uncertain probability distribution $f(\mathbf{v})$ residing in $\mathcal{P}$. 

While the above two formulations have shown to be effective in reducing the conservatism of the optimal solution, in order to implement them, an approximation of the probability distribution $f(\mathbf{v})$ is required, which generally is non-trivial to obtain in practice, especially for high dimensional uncertainty. Besides, an accurate approximation of $f(\mathbf{v})$ generally implies considering a large number of uncertainty samples, which could then induce a large number of constraints in the resulting optimization problems. These issues might limit the applicability of the above formulations. 

Motivated by the objective functions in \eqref{eq:Q_dro} and the property that the uncertainty set consists of multiple subsets, we propose a new objective function
\begin{equation}\label{eq:two_stage_dro_obj}
        \min_{\bf x}\left\{\mathbf{c^\rmt x} + \max_{\mathbf{p}\in\mathcal{P}}\mathbb{E}_{\mathbf{p}}\left[\max_{\mathbf{v}\in\mathcal{V}_k}\min_{\mathbf{y}}\mathbf{b^\rmt y}\right]\right\}
\end{equation}
where $\mathbf{p} = [\mathbf{p}_1,\mathbf{p}_2\cdots,\mathbf{p}_K]^\rmt$ is the probability vector with $\mathbf{p}_k:=\mathbb{P}(\mathbf{v}\in\mathcal{V}_k)$, $\mathcal{P}$ is the ambiguity set defining the admissible set of $\mathbf{p}$. 

\remark{In contrast to the conventional RO formulation in \eqref{eq:ro}, where the worst-case performance among all possible uncertainty scenarios are optimized, the proposed formulation \eqref{eq:two_stage_dro_obj} optimizes the distributionally robust solution w.r.to the expected worst-case performance $\max_{\mathbf{v}\in\mathcal{V}_k}\min_{\mathbf{y}}\mathbf{b}^\rmt\mathbf{y}$ over all uncertainty subsets. In \eqref{eq:two_stage_dro_obj}, the worst-case performance in each uncertainty subset $\mathcal{V}_k$ is weighted based on the probability $\mathbb{P}(\mathbf{v}\in\mathcal{V}_k)$, and hence a less conservative solution is expected with our proposed formulation than the conventional RO formulation. 
}

\remark{In contrast to the SO and DRO formulations in \eqref{eq:obj_options}, which entail the information of the joint probability distribution $f(\mathbf{v})$ of the continuous uncertainty $\mathbf{v}$, our proposed formulation \eqref{eq:two_stage_dro_obj} only considers the discrete probability distribution $\mathbb{P}(\mathbf{v}\in\mathcal{V}_k)$. 
Compared with $f(\mathbf{v})$, $\mathbb{P}(\mathbf{v}\in\mathcal{V}_k)$ is much easier to estimate with high accuracy since $f(\mathbf{v})$ is the continuous probability distribution of random variable $\mathbf{v}$ that might be high dimensional. In contrast, the probability vector $\mathbf{p}$ is the distribution of the discrete random variable $\mathbb{I}(\mathbf{v}\in\mathcal{V}_k)$, where $\mathbb{I}(\cdot)$ is the indicator function. An approximation of $\mathbf{p}$ can be readily computed by counting the frequency of $\mathbf{v}\in\mathcal{V}_k$ for all uncertainty data samples, and is a multinomial distribution when uncertainty subsets $\{\mathcal{V}_k,\ k=1,\cdots,K\}$ are disjoint. Furthermore, the uncertainty is still described via an uncertainty set in \eqref{eq:two_stage_dro_obj}, and uncertainty samples are not explicitly considered, which are beneficial for maintaining the computational efficiency as in conventional RO formulation. On the contrary, the resulting optimization problems for the existing SO and DRO formulations in the format of \eqref{eq:obj_options}, e.g., the scenario-based SO formulation \citep{calafiore2006scenario} and Wasserstein-based DRO formulation \citep{mohajerin2018data}, need to consider all uncertainty samples and are generally computationally demanding in the presence of large numbers of uncertainty scenarios.
}

In this work, the ambiguity set $\mathcal{P}$ is defined based on Kullback Leibler (KL) divergence as
\begin{equation}\label{eq:ambiguity}
    \mathcal{P}:=\left\{\mathbf{p}\mid \mathbf{p} \geq 0,\ \mathbf{e}^\rmt\mathbf{p} = 1,\ KL(\mathbf{\bar p},\mathbf{p})\leq \rho\right\},
\end{equation}
where ${KL}(\mathbf{\bar p},\mathbf{p}):= \sum_{k}\mathbf{\bar p}_k\log (\frac{\mathbf{\bar p}_k}{\mathbf{p}_k})$ is the KL divergence function measuring the similarity of two probability vectors $\mathbf{\bar p}$ and $\mathbf{{p}}$ with $\bar{\bf p}$ as the approximated probability vector extracted from historical data, $\rho \geq 0$ is a user-defined parameter determining the size of the ambiguity set. 

Considering the objective function \eqref{eq:two_stage_dro_obj} as well as the ambiguity set definition \eqref{eq:ambiguity} and applying epigraphical reformulation, it gives to the following optimization problem
\begin{subequations}\label{eq:two_stage_dro}
    \begin{align}
        \min_{\mathbf{x},\eta} \ & \mathbf{c^\rmt x} + \eta  \\
        \text{s.t. } &  \max_{\mathbf{p}\in\mathcal{P}}\bigg\{\mathbb{E}_{\mathbf{p}}\left[\max_{\mathbf{v}\in\mathcal{V}_k}\min_{\mathbf{y}}\mathbf{b}^\rmt \mathbf{y}\right]\bigg\} \leq \eta, \label{eq:max_expect}\\
        &     \mathbf{p} \geq 0,\ \mathbf{e}^\rmt\mathbf{p} = 1,\ KL(\mathbf{\bar p},\mathbf{p})\leq \rho,\\
        & \mathbf{Ax \leq  q},\\
        & \mathbf{Tx + Wy + Mv \leq h}.
    \end{align}
\end{subequations}

\remark{For the proposed formulation in \eqref{eq:two_stage_dro}, the conservatism of the optimal solution can be influenced by the size of the ambiguity set $\mathcal{P}$ and the number of uncertainty subsets $K$. Namely, the conservatism can be adjusted both in constructing the data-driven uncertainty set (the value of $K$) and in designing the ambiguity set (the value of $\rho$).
Roughly speaking, a larger size of $\mathcal{P}$ contributes to more conservatism while a larger $K$ leads to less conservatism. In extreme cases, given a fixed $K$, a sufficiently large size of $\mathcal{P}$ leads to optimizing the worst-case performance as with the RO formulation in \eqref{eq:ro} since the ambiguity set will include the probability vector $\mathbf{p}$ assigning all probability mass on the uncertainty set containing the worst uncertainty sample among all possible uncertainties. Likewise, if the uncertainty set has no subset, the proposed formulation with any feasible size of $\mathcal{P}$ also reduces to worst-case performance optimization as in \eqref{eq:ro}.
}

\subsection{Computational Method Design}
The proposed formulation \eqref{eq:two_stage_dro} contains quadruple-level optimization problems and is computationally intractable to numerical solvers. This section will present a computational method based on the CCG algorithm to efficiently solve the optimization problem \eqref{eq:two_stage_dro}.

\assumption{For the robust optimization problem \eqref{eq:two_stage_dro}, the uncertainty set $\mathcal{V}$ is a union of $K$ nonempty and bounded subsets $\mathcal{V}_k$ defined as $\mathcal{V}_k:=\{\mathbf{v}\mid \mathbf{D}_k\mathbf{v} \leq \mathbf{d}_k\}$.}

\theorem{\textit{Assuming Assumption 2 holds, the optimization problem \eqref{eq:two_stage_dro} can be solved via Algorithm \ref{alg:ccg_dro} by iteratively solving a master problem {\rm \textbf{MP}}$_{\text{\rm DRO}}$ in \eqref{eq:mp} and $K$ subproblems {\rm\textbf{SP}}$_{\text{\rm DRO}}^k$ in \eqref{eq:sp} within finite number of iterations. }} 
\begin{subequations}\label{eq:mp}
    \begin{align}
        \textbf{MP}_{\text{DRO}}:& \notag\\
        \min_{\mathbf{x},\eta,\nu,\mu,\mathbf{y}_{k,i}}\ & \mathbf{c^\rmt x} + \eta \\
        \text{s.t. } &  \nu \sum\nolimits_{k=1}^K \bar{\mathbf{p}}_k \exp\bigg(\frac{\mathbf{b^\rmt y}_{k,i} - \mu}{\nu} -1\bigg) \notag\\
        & \quad +\mu + \rho \nu \leq \eta, \label{eq:mp_exp}\\
        & \mathbf{Ax \leq q}, \\
        & \mathbf{Tx + W}\mathbf{y}_{k,i} + \mathbf{Mv}^*_{k,i} \leq \mathbf{h},\\
        & \nu \geq 0,\\
        & i = 1,\cdots, r,\ k = 1,\cdots,K.
    \end{align}
\end{subequations}

\textbf{SP}$_{\text{DRO}}^k$:
\begin{subequations} \label{eq:sp}
    \begin{align}
        \max_{\mathbf{v}}\min_{\mathbf{y}}\ & \mathbf{b^\rmt y} \\
        \text{s.t. } & \mathbf{Tx}^* + \mathbf{Wy + Mv \leq h}, \\
         & \mathbf{D}_k\mathbf{v}\leq \mathbf{d}_k.
    \end{align}
\end{subequations}
\textit{Proof}. The proof is provided in the Appendix.\hfill $\square$

\remark{It should be noted that the master problem \textbf{MP}$_{\text{DRO}}$ \eqref{eq:mp} is nonlinear due to the constraint \eqref{eq:mp_exp} even if the original constraints and objective function in \eqref{eq:two_stage_dro} are linear. Besides, it can be seen from Algorithm \ref{alg:ccg_dro} that the number of these nonlinear constraints will increase with the iteration of the algorithm, which could degrade the computational efficiency of our proposed scheme. One solution to alleviate this issue is applying the following reformulation for constraint \eqref{eq:mp_exp}:
\begin{subequations}\label{eq:mp_exp_reform}
    \begin{align}
        &\mu + \rho \nu + \nu\sum_k \bar{\mathbf{p}}_k \exp\left(\frac{\phi_k - \mu}{\nu} - 1\right) \leq \eta, \label{eq:17a}\\
        & \mathbf{b}^\rmt \mathbf{y}_{k,i} \leq \phi_k,\ i = 1,\cdots, r,\ k = 1,\cdots,K,\label{eq:mp_exp_reform1}
    \end{align}
\end{subequations}
where $\phi_k$ are auxiliary decision variables. The above reformulation is valid because the left-hand side of \eqref{eq:mp_exp} is strictly increasing in $\mathbf{b}^\rmt \mathbf{y}_{k,i}$. By replacing \eqref{eq:mp_exp} via \eqref{eq:mp_exp_reform}, only linear constraints are added in each iteration, and the number of nonlinear constraints \eqref{eq:17a} is fixed regardless of the algorithm iteration, which is beneficial in improving the computational efficiency of Algorithm \ref{alg:ccg_dro}. 
}
\begin{algorithm}[htb]
\caption{\textit{column-and-constraint generation} algorithm for solving \eqref{eq:two_stage_dro}.}
\label{alg:ccg_dro}
    \begin{algorithmic}[1]
        \Statex \textbf{Input}: suboptimality gap $\epsilon$
        \Statex \textbf{Output}: optimal decision variable $\mathbf{x}^*$ and objective function value $\mathbf{c}^\rmt\mathbf{x}^* + \eta^*$
        \State Set $LB = -\infty$, $UB = \infty$, $r = 0$ 
        \While{$|UB - LB| > \epsilon$} 
        \State Solve \textbf{MP}$_{\text{DRO}}$ \eqref{eq:mp} to obtain solutions $(\mathbf{x}^*, \eta^*, \mu^*, \nu^*$) and update $LB = \mathbf{c^\rmt} \mathbf{x}^*+\eta^*$
        \State Solve $K$ subproblems \textbf{SP}$_{\text{DRO}}^k$ \eqref{eq:sp} to obtain solutions $\{\mathbf{v}^*_k,\mathbf{y}_{k}^*\}_{k=1}^{K}$, and update $UB$ as 
        \begin{equation*}
        \begin{aligned}
            UB = \min & \bigg\{UB,\ \mathbf{c^\rmt} \mathbf{x}^* + \mu^* + \rho \nu^* + \\&\nu^*\sum\nolimits_{k=1}^{K}\bar{\mathbf{p}}_k\exp\left(\frac{\mathbf{b}^\rmt \mathbf{y}_k^{*} - \mu^*}{\nu^*}-1\right)\bigg\}
            \end{aligned}
        \end{equation*}
        \State {Create decision variables $\{\mathbf{y}_{k,r}\}_{k=1}^{K}$, set parameters $\mathbf{v}^*_{k,r} = \mathbf{v}_k^*\ (k=1,\cdots,K)$, and add the following constraints to \textbf{MP}$_{\text{DRO}}$ in \eqref{eq:mp}}
        \begin{align*}
        \begin{cases}
            &\mu+\rho \nu + \nu\sum\nolimits_{k=1}^K\bar{\mathbf{p}}_k\exp\left(\frac{\mathbf{b}^\rmt \mathbf{y}_{k,r} - \mu}{\nu}-1\right) \leq \eta ,\\
            &\mathbf{Tx + Wy}_{k,r} +\mathbf{Mv}^*_{k,r} \leq \mathbf{h},\ k = 1,\cdots,K.
            \end{cases}
        \end{align*}
        \State $r \leftarrow r+1$
        \EndWhile
        \State \textbf{Return}: $\mathbf{x}^*$ and $\mathbf{c^\rmt}\mathbf{x}^* + \eta^*$
    \end{algorithmic}
\end{algorithm}

\remark{It is worth noting that, while the constraints \eqref{eq:mp_exp} and \eqref{eq:17a} are nonlinear, they are convex, and \textbf{MP}$_{\text{DRO}}$ becomes a convex optimization when the decision variable $\mathbf{x}$ is continuous. The convexity of the constraints \eqref{eq:mp_exp} and \eqref{eq:17a} is proved in the Appendix. Several off-the-shelf solvers, such as {\tt Gurobi}, {\tt Ipopt} and {\tt MadNLP}, can deal with this type of nonlinearity. The subproblem \textbf{SP}$_{\text{DRO}}$ \eqref{eq:sp} is a bi-level linear optimization problem that is not numerically tractable to solvers. However, by applying strong duality or KKT-based reformulations and big-M approach, the inner-level optimization problem in \eqref{eq:sp} can be eliminated, and this bilevel optimization problem can be reformulated as a mixed-integer linear programming problem, see \cite{zeng2013solving,zhao2012exact} for more details. 
}

\remark{Compared to conventional RO formulations, both the proposed formulation \eqref{eq:two_stage_dro} and existing DRO formulations, e.g., Wasserstein-based DRO formulation \citep{mohajerin2018data}, result in increased computational burden, though for different reasons. In our proposed formulation, this computational burden arises primarily from the nonlinear constraint \eqref{eq:mp_exp}, particularly when the first-stage decision variables $\mathbf{x}$ include binary components, making the master problem \textbf{MP}$_{\text{DRO}}$ \eqref{eq:mp} a mixed-integer nonlinear optimization problem. As for the existing DRO formulation, we take the Wasserstein-based DRO \citep{mohajerin2018data} as an example. When considering the fully adaptive recourse decision rule and large numbers of uncertainty samples, while the resulting optimization problem still retains the linearity of the original deterministic formulation, the increased computational burden is mainly caused by solving numerous subproblems in every iteration of the CCG algorithm since each uncertainty scenario will incur an optimization problem, see Algorithm 1 in the supplementary material \cite{yun24}. The curse of dimensionality issue of the Wasserstein-based DRO formulation is also mentioned in \cite{wang2020wasserstein}. }

\section{Simulation Results}\label{sec:simulation}
This section presents three case studies to illustrate the effectiveness of the proposed schemes in this paper. \textit{Case Study 1} considers robust predictive control of building climate to demonstrate the effectiveness of the proposed approach in Section \ref{sec:predictive_control}. \textit{Case Study 2} and \textit{Case Study 3} showcase the approach designed in Section \ref{sec:stochastic_dro} with robust location transportation planning and chemical process network planning problems, respectively.

All simulations are implemented on an Intel Xeon W-2223 CPU at 3.6GHz with 16GB RAM. Optimization problems are modelled via Python package {\tt gurobipy} and solved via {\tt Gurobi 11.0} \citep{gurobi}. The values of the parameters used in our case studies are provided in the supplementary material \cite{yun24}.

\subsection{Case Study 1: Robust Predictive Control of Building Climate}
This case study considers robust predictive control of building climate. Building systems suffer from weather uncertainties, such as prediction errors of ambient temperature, solar irradiation, etc. Properly considering these uncertainties can improve indoor climate comfort. 

A building climate predictive control problem in its deterministic form can be formulated as \citep{shang2019data}
\begin{eqnarray*}
    \min_{\mathbf{v}_t}\ && \sum_{t=1}^Nl(\mathbf{s}_t,\mathbf{v}_t)\\
    \text{s.t. }&& \mathbf{s}_{t+1} = \bm{\Phi}\mathbf{s}_t + \bm{\Gamma}_u \mathbf{u}_t + \bm{\Gamma}_w\mathbf{w}_t + \bm{\Gamma}_v\mathbf{v}_t,\\
    && \mathbf{\underline{s}}_t \leq \mathbf{s}_t \leq \bar{\mathbf{s}}_t,\quad \underline{\mathbf{u}}_t \leq \mathbf{u}_t \leq \bar{\mathbf{u}}_t, \\
    && \forall \mathbf{v}_t \in\mathcal{V}_t,\ t = 1,\cdots,N
\end{eqnarray*}
where $l(\cdot,\cdot)$ is the stage cost function, $\mathbf{s}_t$ is the system states consisting of indoor temperature, roof temperature, wall temperature and floor temperature; $\mathbf{u}_t$ denotes the heating power, $\mathbf{w}_t$ is the predicted ambient conditions, $\mathbf{v}_t$ is the prediction error of ambient temperature, $\mathcal{V}_t$ denotes the uncertainty set, $N$ is the length of prediction horizon, $\underline{\mathbf{s}}_t$/$\underline{\mathbf{u}}_t$ and $\bar{\mathbf{s}}_t$/$\bar{\mathbf{u}}_t$ are lower bound and upper bound of system states/control inputs, respectively. The stage cost function is defined as $l(\mathbf{s}_t,\mathbf{u}_t):= \mathbf{u}_t$ to minimize energy usage, system state constraints are defined to keep the indoor temperature above $21^\circ$C during 7:00 - 18:00 and above 15$^\circ$C during remaining hours. Heating power constraints are $0 \leq \mathbf{u}_t \leq 150$. The values of system matrices $(\bm{\Phi},\bm{\Gamma}_u,\bm{\Gamma_w},\bm{\Gamma}_v)$ are adopted from \cite{shang2017data}.

In the simulation, uncertainties $\mathbf{v}_t$ $(t=1,\cdots,N)$ are assumed to be I.I.D., and the uncertainty set $\mathcal{V}_t$ for $\mathbf{v}_t$ is a union of 2 subsets. As a result, the total number of subsets for the stacked uncertainty $\mathbf{v} = [\mathbf{v}_1^\rmt,\cdots,\mathbf{v}_N^\rmt]^\rmt$ is $2^N$. 
Two schemes are considered in this case study for solving the building climate control problem. One is the conventional RO formulation with the explicit description of each uncertainty subset \eqref{eq:uncerset_exp} and is solved via the CCG-based algorithm in \cite{ning2017data}. Another one is our proposed formulation \eqref{eq:RO_reform} and is solved via Algorithm \ref{alg:ccg_worst_case}. To fully demonstrate the computational efficiency of the proposed scheme, different values of the prediction horizon $N$ are tested. 
\begin{figure}[htb]
    \centering
    \begin{minipage}{\linewidth}
    \centering
    \subfloat[computational time.]{\includegraphics[width=0.9\linewidth]{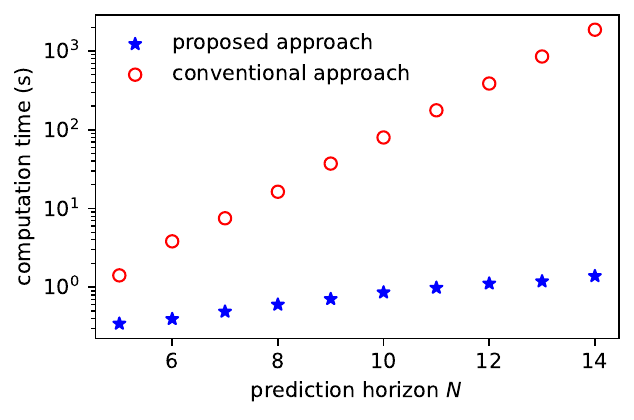}}
    \end{minipage}
    \begin{minipage}{\linewidth}
    \centering
    \subfloat[objective values]{\includegraphics[width=0.9\linewidth]{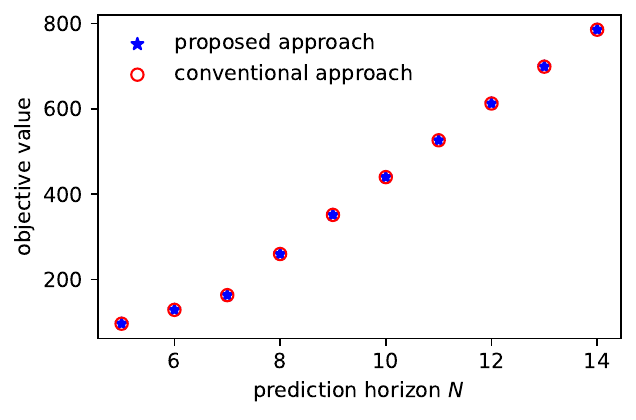}}
    \end{minipage}
    \caption{Simulation results of \textit{Case Study 1} with the proposed approach and the conventional approach.}
    \label{fig:case1}
\end{figure}
Simulation results are shown in Figure \ref{fig:case1}. It can be seen that, with the increase of the prediction horizon $N$, the computational time of applying the conventional approach increases exponentially and is much larger than that of our proposed formulation since $2^N$ numbers of subproblems have to be solved in each algorithm iteration. In contrast, the computation time with our proposed formulation \eqref{eq:RO_reform} and Algorithm \ref{alg:ccg_worst_case} only increases linearly with the prediction horizon $N$. Besides, Figure \ref{fig:case1}(b) shows that both two approaches give the same objective function, which indicates that the proposed formulation does not sacrifice optimality while remarkably improving the computational efficiency.
\begin{table*}[htb]
    \centering
        \caption{Computational results of different RO formulations for \textit{Case Study 2}.}
    \label{tab:location_transportation}
    \resizebox{0.75\linewidth}{!}{\begin{tabular}{l|c|c|c|c}\toprule\hline
    & \multirow{2}{*}{\makecell[c]{conventional RO\\formulation\\(Scheme 1)}} &\multirow{2}{*}{\makecell[c]{existing DRO formulation\\with Wasserstein metric\\(Scheme 2) }} &  \multicolumn{2}{c}{the proposed RO formulation \eqref{eq:two_stage_dro}}\\\cline{4-5}
        & &  & \makecell{with \eqref{eq:mp_exp}\\(Scheme 3)} & \makecell{with \eqref{eq:mp_exp_reform}\\(Scheme 4)} \\\hline
        Optimal Objective & 36632& 35238  & 35419 & 35419 \\
        Decision Variable $\mathbf{x}$ & [1,0,1,274,0,570] & [1,0,1,324,0,520] & [1,0,1,522,0,322] & [1,0,1,522,0,322] \\
        CPU Time (s) & 0.40 & 425.80 & 11.14 & 8.07 \\
        Iterations & 2& 3 & 2& 2 \\\hline
    \end{tabular}}
\end{table*}
\subsection{Case Study 2: Robust Location Transporation Planning}
In this subsection, the effectiveness of the proposed formulation \eqref{eq:two_stage_dro} and the corresponding Algorithm \ref{alg:ccg_dro} is validated via a robust location transportation planning problem, which is also considered as a benchmark problem in \cite{zeng2013solving}. The deterministic formulation of this problem is
\begin{subequations}
    \begin{align}
        \min\ & 400x_1 + 414x_2+326x_3 + 18x_4 + 25 x_5 + 20x_6 \notag\\
        & + 22y_{11} + 33y_{12} + 24y_{13} + 33y_{21} + 23y_{22} + 30y_{23}\notag\\
        & 20y_{31} + 25y_{32} + 27y_{33} \\
        \text{s.t. } &  x_{i+3} \leq 800x_i,\quad i = 1,2,3,\\
        & \sum_{j} y_{ij} \leq x_{i+3},\quad \forall i = 1,2,3,\\
        & \sum_iy_{ij}\geq d_j + 40*v_j \\
        & x_i\in\mathbb{B},\quad  x_{i+3}\leq 0, \quad i = 1,2,3,\\
        & y_{ij} \geq 0,\quad \forall i = 1,2,3,\quad  j = 1,2,3
    \end{align}
\end{subequations}
where binary first-stage decision variables $(x_1,x_2,x_3)$ determine the location of the facilities; continuous first-stage decision variables $(x_4,x_5,x_6)$ denote the facility capacities; recourse decision variables $y_{ij}\ (i=1,2,3,\ j=1,2,3)$ are transportations, $[d_1,d_2,d_3] = [206,274,220]$ are basic demands, and $v_j\ (j=1,2,3)$ are scaled demand uncertainties. 

In our simulation, four schemes are considered: 
\begin{itemize}
    \item Scheme 1: the conventional RO formulation \eqref{eq:ro}.
    \item Scheme 2: the Wasserstein-based DRO formulation \citep{mohajerin2018data}. 
    \item Scheme 3: the proposed formulation \eqref{eq:two_stage_dro} and Algorithm \ref{alg:ccg_dro} with the constraints \eqref{eq:mp_exp}.
    \item Scheme 4: the proposed formulation \eqref{eq:two_stage_dro} and Algorithm \ref{alg:ccg_dro} with the reformulated constraints \eqref{eq:mp_exp_reform}.
\end{itemize}

The uncertainty set $\mathcal{V}$ for $[v_1,v_2,v_3]^\rmt$ is supposed to have 4 polyhedral subsets. For our proposed schemes (Scheme 3 \& 4), the nominal probability distribution for defining the ambiguity set $\mathcal{P}$ is $\bar{\bf p} = [0.7, 0.1, 0.1, 0.1]$. The upper bound of KL divergence for defining the ambiguity set $\mathcal{P}$ in \eqref{eq:ambiguity} is set as $\rho = 0.5$. For Scheme 2, since its implementation entails considering specific uncertainty scenarios, 1000 uncertainty samples in $\mathcal{V}$ are randomly generated, which is a reasonable choice for the uncertainty with dimension 3 \citep{van1999models,scott2015multivariate}. Since we mainly focus on demonstrating the computational efficiency of the DRO formulation, the Wasserstein distance $\varepsilon$, which does not affect its computational efficiency, is set as $\varepsilon=1$ to give a comparable performance as with the other schemes.

Simulation results are summarized in Table \ref{tab:location_transportation}. It can be seen that these different RO formulations derive distinct first-stage decision variables with non-trivial differences. Compared with the conventional RO formulation (Scheme 1) for optimizing the worst-case performance, the remaining formulations give less conservative solutions, i.e., smaller optimal objective values, at the price of increased computational burden. As discussed in Remark 8, the increased computational burden for the DRO formulation (Scheme 2) is due to its curse of dimensionality when considering the fully adaptive decision rule and a large number of uncertainty samples \citep{mohajerin2018data,wang2020wasserstein}, while for our proposed schemes (Scheme 3 \& 4), the increased computational burden is caused by the nonlinear constraint \eqref{eq:mp_exp}, especially when the first-stage decision variable $\mathbf{x}$ contains integer ingredients. It can be seen that, in comparison with the DRO formulation (Scheme 2), our proposed formulations (Scheme 3 \& 4) are much more computationally efficient. Furthermore, the computational time of Scheme 4 is less than that of Scheme 3, which indicates the efficacy of the proposed reformulation \eqref{eq:mp_exp_reform} for improving the computational efficiency of Algorithm \ref{alg:ccg_dro}. 

\subsection{Case Study 3: Chemical Process Network Planning}
To further illustrate the viability and effectiveness of the proposed formulation \eqref{eq:two_stage_dro} and Algorithm \ref{alg:ccg_dro}, a chemical process network planning (CPNP) problem is investigated. The CPNP problem is a typical engineering problem that fits into RO settings and has been considered as benchmark problems in several existing literature, see \cite{you2011stochastic,shang2017data,ning2018data}. 

\begin{figure}[tb]
    \centering
    \includegraphics[width=\linewidth]{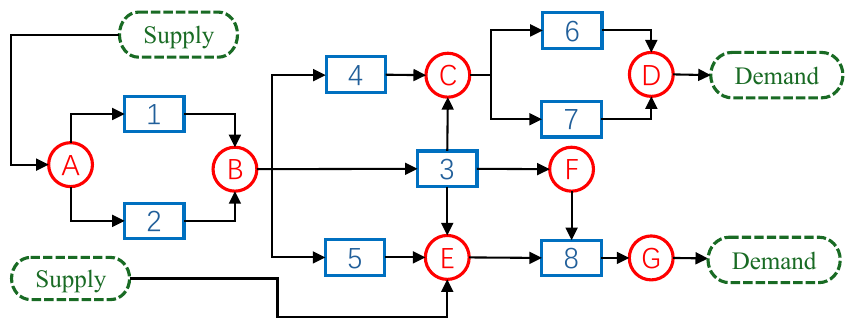}
    \caption{The chemical process network for \textit{Case Study 3}.}
    \label{fig:chem_process}
\end{figure}

A chemical process consists of raw materials, intermediate chemicals, final products, and multiple interconnected processes. The design objective of CPNP is to maximize the net present value (NPV) of the entire network while respecting system constraints for all possible uncertainties. Our case study considers a chemical process network consisting of 8 processes and 7 chemicals. The network sketch is shown in Figure \ref{fig:chem_process}, where chemicals are denoted as red circles ($A,B,\cdots,G$), processes are blue rectangles ($1,2,\cdots,8$), and process flows are indicated as arrows. Among all chemicals, $(A,E)$ are raw materials, and $(D,G)$ are products. 

The RO problem for CPNP within our proposed design framework \eqref{eq:two_stage_dro} can be formulated as
\begin{subequations}\label{eq:chem_process}
    \begin{align}\hspace{-5pt}
        \max_{\substack{QE_{i,t},Y_{i,t}}} &\sum_{i}\sum_{t}(-\alpha_{i,t}  QE_{i,t} - \beta_{i,t} Y_{i,t}) + \notag\\        &\min_{\mathbf{p}\in\mathcal{P}}\mathbb{E}_{\mathbf{p}}\bigg[\min_{\mathbf{v}\in\mathcal{V}_k} \max_{\substack{P_{j,t},QE_{i,t}\\ S_{j,t},W_{i,t} }}\bigg(-\sum_{i}\sum_t \gamma_{i,t}W_{i,t}\notag\\
        &-  \sum_{j}\sum_{t}\varphi_{j,t} P_{j,t} + \sum_{j}\sum_{t}\tau_{j,t}S_{j,t}\bigg)\bigg] \label{eq:chem_obj}\\
        \text{s.t. } & qe_{i,t}^L\cdot Y_{i,t} \leq QE_{i,t} \leq qe_{i,t}^U\cdot Y_{i,t},\ \forall i\in\mathcal{I},\ \forall t\in\mathcal{T}\label{eq:chem_c1}\\
        & Q_{i,t} = Q_{i,t-1} + QE_{i,t},\ \forall i\in\mathcal{I},\ \forall t \in\mathcal{T}\label{eq:chem_c2}\\
        & \sum_{t}Y_{i,t} \leq ce_i,\ \forall i\in\mathcal{I}\label{eq:chem_c3}\\
        & \sum_{i}\alpha_{i,t}\cdot QE_{i,t} + \beta_{i,t}Y_{i,t} \leq cb_t,\ \forall t \in\mathcal{T}\label{eq:chem_c4}\\
        & W_{i,t} \leq Q_{i,t}, \ \forall i \in\mathcal{I},\ \forall t \in\mathcal{T}\label{eq:chem_c5}\\
        & P_{j,t} - \sum_i\kappa_{i,j}W_{i,t} - S_{j,t} = 0,\ \forall j \in\mathcal{J},\ \forall t\in\mathcal{T}\label{eq:chem_c6}\\
        & P_{j,t} \leq su_{j,t},\ S_{j,t} \leq du_{j,t},\forall j \in\mathcal{J},\ \forall t \in\mathcal{T}\label{eq:chem_c7}\\
        & QE_{i,t}, Q_{i,t}, P_{j,t}, W_{i,t}, S_{j,t} \geq 0,Y_{i,t} \in\{0,1\}, \notag\\
        &\quad \forall i\in\mathcal{I},\ \forall j \in\mathcal{J},\ \forall t\in\mathcal{T}, \label{eq:chem_c8}\\
        & \mathbf{v} = \{du_{i,t},su_{i,t}\} \in\mathcal{V},\ \forall i \in \mathcal{I},\ \forall t\in\mathcal{T}
    \end{align}
\end{subequations}\noindent where $\mathcal{I} = \{1,\cdots,I\}$, $\mathcal{J}= \{1,\cdots,J\}$, and $\mathcal{T} = \{1,\cdots,T\}$ with $I,J$ and $T$ as the total numbers of processes, chemicals and planning periods, respectively. All notations in the above equations are explained in Table \ref{tab:chem_notations}. The objective function \eqref{eq:chem_obj} maximizes the NPV consisting of investment cost, operation cost, purchase cost of raw chemicals, and sale profit; constraint \eqref{eq:chem_c1} specifies the upper and lower bounds of capacity expansion for all processes and time periods; constraint \eqref{eq:chem_c2} updates the total available capacity of each process; constraint \eqref{eq:chem_c3} limits the largest process expansion times; constraint \eqref{eq:chem_c4} ensures the process expansion costs are within available budgets; inequality \eqref{eq:chem_c5} limits the production level of each process within its total capacity; equality \eqref{eq:chem_c6} models the mass balance of all chemicals; constraints \eqref{eq:chem_c7} ensures that the amounts of purchased and sold chemicals are limited by the available market supply and demand, respectively; constraints \eqref{eq:chem_c8} indicates all non-negative continuous decision variables and binary decision variables. For more detailed explanations of chemical process networks, please refer to \cite{you2011stochastic,ning2018data}. 
\begin{table}[htb]
\caption{Notations in \eqref{eq:chem_process}}
    \centering\resizebox{\linewidth}{!}{  
    \begin{tabular}{ll ll}\toprule\hline
        \multicolumn{4}{c}{decision variables} \\\hline
         $QE_{i,t}$ & \makecell[l]{amount of capacity\\ expansion} & $Y_{i,t}$ & \makecell[l]{binary decision of\\capacty expansion} \\
         $P_{j,t}$ & chemical purchase amount & $S_{j,t}$ & chemical sale amount \\
         $W_{i,t}$ & operating level & $Q_{i,t}$ & total process capacity \\
         $su_{i,t}$ & market supply limit & $du_{i,t}$ & market demand limit\\\hline
         \multicolumn{4}{c}{parameters} \\\hline
         $\alpha_{i,t}$& variable investment cost & $\beta_{i,t}$ & fixed investment cost \\
         $\gamma_{j,t}$& operating cost & $\varphi_{j,t}$ & purchase cost\\
         $\tau_{j,t}$ & sale price & $qe_{i,t}^L$ & \makecell[l]{capacity expansion\\lower bound}\\
         $qe_{i,t}^U$ & \makecell[l]{capacity expansion\\upper bound} & $ce_i$ & expansion number limit \\
         $cb_t$ & expansion cost budget & $\kappa_{i,j}$ & mass balance coefficient\\
         $i$ & index of the $i$-th process & $j$ & index of the $j$-th chemical \\
         $t$ & index of the $t$-th time period \\\hline
    \end{tabular}}
    \label{tab:chem_notations}
\end{table}

\begin{table*}[htb]
    \centering
        \caption{Computational results of different RO formulations for \textit{Case Study 3}.}
    \label{tab:cpnp_results}
    \resizebox{0.7\linewidth}{!}{
    \begin{tabular}{l|c|c|c|c}\toprule\hline
    & \multirow{2}{*}{\makecell[c]{conventional RO\\formulation\\(Scheme 1)}} & \multirow{2}{*}{\makecell[c]{existing DRO formulation\\with Wasserstein metric\\(Scheme 2)}} & \multicolumn{2}{c}{the proposed RO formulation \eqref{eq:two_stage_dro}}\\\cline{4-5}
        & & & \makecell{with \eqref{eq:mp_exp}\\(Scheme 3)} & \makecell{with \eqref{eq:mp_exp_reform}\\(Scheme 4)} \\\hline
        Max. NPV (\$MM) & 221 & 284 &288 &288  \\
        CPU Time (s) & 17.63 & 8758.07 &195.86 &183.17 \\
        Iterations & 2& 5 & 2& 2 \\\hline
    \end{tabular}} 
\end{table*}

The uncertainty variables considered in our design are the market supply limit and demand limit $(su_{i,t},du_{i,t})$ for all raw materials and products over the planning period. Expansion decisions $(QE_{i,t},Y_{i,t})$ are first-stage decision variables, and other remaining variables $(P_{j,t},S_{j,t},Q_{i,t},W_{i,t})$ are recourse decision variables. Our case study considers a $5-$year planning period with each planning period as 1 year. The uncertainties $(su_{j,t},du_{j,t})$ are assumed to be independent and reside in 4 uncertainty sets.

As in \textit{Case Study 2}, we consider the same four schemes. For the DRO formulation in Scheme 2, 1000 uncertainty scenarios are randomly generated, and the Wasserstein distance $\varepsilon$ defining the size of the distribution ambiguity set is $\varepsilon = 300$ to give a comparable objective value with the other schemes. It should be noted that the number of uncertainty samples considered is far smaller than needed for a proper approximation of the uncertainty probability distribution, given the 20-dimensional uncertainty space \citep{van1999models,scott2015multivariate}. However, as shown in our simulation results, even this limited number of samples can be computationally demanding for solving the corresponding DRO problem. For our proposed schemes (Scheme 3 \& 4), the parameters defining the ambiguity set $\mathcal{P}$ are $\bar{\mathbf{p}} = [0.5,0.1,0.2,0.2]$ and $\rho = 0.5$.

Simulation results are shown in Table \ref{tab:cpnp_results}, from which it can be concluded that consistent results are obtained as in \textit{Case Study 2}. Compared with the conventional RO formulation (Scheme 1), the remaining formulations (Scheme 2, 3 \& 4) give less conservative solutions, i.e., a larger value of Max.~NPV at the price of increased computational burden (longer CPU time). Compared with the Wasserstein-based DRO formulation (Scheme 2), our proposed formulations (Scheme 3 \& 4) are less computationally demanding. In addition, compared with the original constraints \eqref{eq:mp_exp}, the reformulation in \eqref{eq:mp_exp_reform} is effective in improving the computational efficiency of Algorithm \ref{alg:ccg_dro}. 

\begin{figure}[htb]
    \centering
    \includegraphics[width=\linewidth]{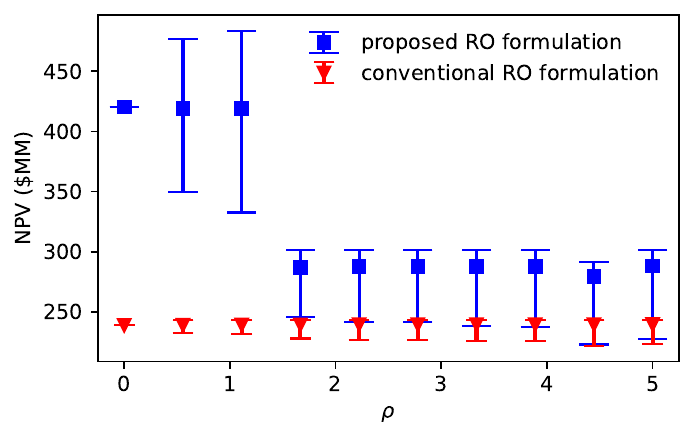}
    \caption{Error bar of the expected NPV with different ambiguity set size $\rho$.}
    \label{fig:errorbar_npv}
\end{figure}

Furthermore, different values of the ambiguity set size $\rho$ are tested to analyze its influence on the conservatism of the optimal solution for our proposed formulation. In the simulation, for each $\rho$, 5000 random trials are implemented to sample the probability distribution $\mathbf{p}$ that resides in the corresponding ambiguity set. For all admissible samples of $\mathbf{p}$, the corresponding objective function values are computed. Figure \ref{fig:errorbar_npv} depicts the error bar (mean value and envelope) of the expected NPV $\mathbf{c}^\rmt\mathbf{x}^* + \sum_{k}\mathbf{p}_k\cdot \mathbf{b}^\rmt \mathbf{y}^*_k$ among all feasible samples of $\mathbf{p}$. It can be seen that, compared with the conventional RO formulation, our proposed approach consistently leads to less conservative solutions for all $\rho$. Namely, the proposed approach gives a higher mean, minimal, and maximal NPV than the conventional RO formulation. Besides, it can be seen from Figure \ref{fig:errorbar_npv} that, with the increase of $\rho$ (the size of the ambiguity set), the minimal NPV of the proposed formulation \eqref{eq:two_stage_dro} converges to that of the conventional RO formulation, which is also consistent with our statement in Remark 5.

\section{Conclusions}\label{sec:conclusion}
This paper investigates the RO problems where uncertainty sets comprise multiple subsets, focusing on two separate questions: 1) \textit{how to address the computational challenge posed by exponentially increasing uncertainty subsets in RO-based predictive control?} and 2) \textit{how to mitigate the conservatism of the robust optimal solutions by leveraging the structure of the uncertainty set?}

To address the first question, we propose a monolithic mixed-integer representation of the uncertainty set. Unlike conventional formulations requiring a separate optimization problem for each subset, our method only solves a single mixed-integer optimization problem to compute the worst-case uncertainty scenario over all subsets. This method is particularly advantageous for RO-based predictive control, where the number of uncertainty subsets could increases exponentially with the prediction horizon.

For the second question, we formulate an innovative objective function exploiting the structure of the multi-subset uncertainty set by combining the existing RO and DRO formulations. The proposed formulation achieves less conservative solutions than conventional RO formulations while showing more computational efficiency than conventional DRO formulations. Besides, a CCG-based algorithm is developed to solve the resulting optimization problem efficiently.

Numerical experiments related to the above research questions are conducted to extensively demonstrate the effectiveness of the proposed schemes. 

\section*{Acknowledgements}
The work was supported by the Brains4Buildings project under the Dutch grant programme for Mission-Driven Research, Development and Innovation (MOOI).

\section*{Appendix}\label{sec:appendix}
{\color{black}
\subsection*{Proof of Theorem 1}
\textit{Proof}: The main difference between the proposed Algorithm \ref{alg:ccg_worst_case} with the conventional CCG algorithm in \cite{zeng2013solving,ning2017data} is the definition of the subproblem \eqref{eq:subp_reform1}. With the mixed-integer representation of the uncertainty set $\mathcal{V}$, each set of feasible binary variables $\{\delta^*_{t,k}\mid t= 1,\cdots,N, k = 1,\cdots,K\}$ uniquely defines a subset of $\mathcal{V}$, and vice versa. As a result, by solving the subproblem \textbf{SP}1 \eqref{eq:subp_reform1} or \textbf{SP}2 \eqref{eq:subp_reform2}, the worst-case uncertainty scenario $\mathbf{v}^*=[\mathbf{v}_1^{*\rmt},\cdots,\mathbf{v}_N^{*\rmt}]^\rmt$ among all subsets together with the corresponding uncertainty subset, defined by $\{\delta^*_{t,k}\mid t= 1,\cdots,N, k = 1,\cdots,K\}$, that contains the uncertainty scenario $\mathbf{v}^*$ will be computed. Based on Assumption 1, it can be readily concluded that there are finite vertices of the uncertainty set $\mathcal{V}$. Then, it follows a similar proof as shown in \cite{zeng2013solving} that Algorithm \eqref{alg:ccg_worst_case} will solve the RO problem \eqref{eq:RO_reform} within finite iterations. This completes the proof.\hfill $\square$

\subsection*{Proof of Theorem 2}
\textit{Proof}: For the constraint \eqref{eq:max_expect}, it can be rewritten as
\begin{equation}\label{eq:expect_sum}
\max_{\mathbf{p}\in\mathcal{P}}\sum_{k=1}^K\mathbf{p}_k\cdot\left(\max_{\mathbf{v}\in\mathcal{V}_k}\mathbf{b}^\rmt\mathbf{y}\right)
\end{equation}
Based on \textit{Lemma 19.1} in \cite{bertsimas2022robust},  the inequality \eqref{eq:expect_sum} is equivalent to finding $\mu$ and $\nu \geq 0$ such that
\begin{equation}\label{eq:20}
        \mu + \rho \nu + \nu\sum_k\bar{\mathbf{p}}_k \exp\big(\frac{\mathcal{C}(\mathbf{x}, \mathcal{V}_k) - \mu}{\nu} -1 \big) \leq \eta 
\end{equation}
where 
\begin{subequations}
    \begin{align}
        \mathcal{C}(\mathbf{x},\mathcal{V}_k) := \max_{\mathbf{v}\in\mathcal{V}_k}\min_{\mathbf{y}}\ & \mathbf{b^\rmt y}\\
        \text{s.t. } & \mathbf{Tx + Wy + Mv \leq h}
    \end{align}
\end{subequations}
Assumption 2 implies that there are finite vertices of the uncertainty set $\mathcal{V}$. Besides, it can be seen that the left-hand side of \eqref{eq:20} is a monotonously increasing function w.r.t. $\mathcal{C}(\mathbf{x},\mathcal{V}_k)$, which is LP w.r.t. $\mathbf{v}$. Consequently, it can be concluded that the optimal $\mathbf{v}$ are taken from the vertices of $\mathcal{V}$. By listing all finite vertices of each uncertainty subset $\mathcal{V}_k$, denoted as $\mathcal{S}_k:=\{\mathbf{v}_{k,1},\cdots,\mathbf{v}_{k,H_k}\}$ with $H_k$ as the number of all vertices, the RO problem \eqref{eq:two_stage_dro} can be rewritten as
\begin{subequations}\label{eq:two_stage_universal}
    \begin{align}
        \min_{\substack{\mathbf{x},\mathbf{y}_{k,i_k}\\\eta,\mu,\nu}} & \ \mathbf{c}^\rmt \mathbf{x} + \eta \\
        \text{s.t. } & \mu + \rho\nu + \nu\sum_{k}\bar{\mathbf{p}}_k\exp\left(\frac{\mathbf{b}^\rmt\mathbf{y}_{k,{i_k}}-\mu}{\nu} - 1\right)\leq \eta,\\
        & \mathbf{Ax \leq q},\\
        & \mathbf{Tx+Wy}_{k,{i_k}} + \mathbf{Mv}_{k,{i_k}} \leq \mathbf{h},\\
        & \forall i_k \in \{1,\cdots,H_k\},\ \forall k \in \{1, \cdots, K\}
    \end{align}
\end{subequations}
where $\mathbf{y}_{k,i_k}$ is the optimal recourse decision variable w.r.t. the uncertainty scenario $\mathbf{v}_{k,i_k}$ $(i=1,\cdots,H_k)$. Namely, 
\begin{equation*}    \mathbf{y}_{k,i_k}=\arg\min_{\mathbf{y}}\mathbf{b^\rmt y} \text{ s.t. } \mathbf{Tx+Wy + Mv}_{k,i_k} \leq \mathbf{h}.
\end{equation*}

Since there are finite extreme uncertainty scenarios, the optimization problem \eqref{eq:two_stage_universal} can be computed by iteratively listing all possible extreme uncertainty scenarios as in the conventional CCG algorithm \citep{zeng2013solving,zhao2012exact}, which leads to Algorithm \ref{alg:ccg_dro}.

In the following, we will prove that Algorithm \ref{alg:ccg_dro} can terminate within finite iterations. Namely, Algorithm \ref{alg:ccg_dro} will either find out all extreme uncertainty scenarios, or terminate when a repeated uncertainty scenario is observed. Assuming at the $r$-th iteration of the algorithm, the extreme uncertainty scenarios $\{\mathbf{v}^*_1,\cdots,\mathbf{v}_K^*\}$ computed by solving $K$ subproblems \eqref{eq:sp} are observed in a previous iteration $t\ (t\leq r-1 )$, then $UB=LB$ and the algorithm terminates. 

Suppose that the optimal decision variables by solving \textbf{MP}$_{\text{DRO}}$ \eqref{eq:mp} at the $r$-th iteration are $(\mathbf{x}^*,\mu^*,\eta^*,\nu^*)$, which further lead to the optimal decision variables of $K$ subproblems \textbf{SP}$_{\text{DRO}}^k$ $\{\mathbf{v}_{1}^*,\cdots,\mathbf{v}_K^*,\mathbf{y}_1^*,\cdots,\mathbf{y}_K^*\}$. It readily gives that 
\begin{align*}
    LB \leq UB \leq &\mathbf{c}^\rmt\mathbf{x}^* + \mu^* +\rho\nu^*+\\
    &\quad \nu^*\sum_{k}\mathbf{\bar{p}}_k\exp\left(\frac{\mathbf{b}^\rmt\mathbf{y}_k^* - \mu^*}{\nu^*} - 1\right)
\end{align*}
Since the worst-case uncertainty scenarios $\{\mathbf{v}_1^*,\cdots,\mathbf{v}^*_K\}$ was observed at a previous iteration, the optimal solution of the master problem \textbf{MP}$_{\text{DRO}}$ at $(r+1)$-th iteration will be the same as with the $r$-th iteration and satisfy
\begin{equation*}
    LB \geq \mathbf{c}^\rmt\mathbf{x}^* + \mu^* + \rho\nu^* + \nu^* \sum_k\mathbf{\bar{p}}_k\exp\left(\frac{\mathbf{b}^\rmt\mathbf{y}_k^* -\mu^* }{\nu^*}-1\right)
\end{equation*}
As a result, it can be concluded that $LB = UB$ at the $(r+1)$-th iteration, and the algorithm terminates. This completes the proof. \hfill $\square$

\subsection*{Proof of the Convexity of Constraint \eqref{eq:mp_exp}}
\textit{Proof}: For brevity of analysis, it is sufficient to prove that the nonlinear function $f(x,y):= x\exp(\frac{y}{x}-1)$  with $x\geq 0$ is convex since the nonlinear constraint \eqref{eq:mp_exp} can be reformulated as
\begin{subequations}
    \begin{align}
        &\mu + \rho\nu + \sum_{k}\mathbf{\bar{p}}_kf(\nu,g_{k,i})\leq \eta,\\
        & g_{k,i} = \mathbf{b}^\rmt \mathbf{y}_{k,i} -\mu.
    \end{align}
\end{subequations}

For the nonlinear function $f(x,y):= x\exp{(\frac{y}{x}-1)}$, its Hessian matrix is
\begin{equation}
    \nabla^2 f = 
    \begin{bmatrix}
    \frac{y^2}{x^3}\exp(\frac{y}{x}-1) & -\frac{y}{x^2}\exp(\frac{y}{x}-1) \\
    -\frac{y}{x^2}\exp(\frac{y}{x}-1) & \frac{1}{x}\exp(\frac{y}{x}-1)
    \end{bmatrix}.
\end{equation}
It can be easily verified that $\nabla^2 f$ is positive semidefinite when $x\geq 0$. As a result, it follows from the \textit{second-order condition} of convexity \citep{boyd2004convex} that the nonlinear function $f(x,y)$ is convex, and hence constraint \eqref{eq:mp_exp} is convex. This completes the proof. \hfill $\square$

\bibliographystyle{elsarticle-harv}
\bibliography{ref}

\end{document}